\documentclass{statsoc}
\usepackage{amsmath}
\usepackage{latexsym}
\usepackage{amssymb}
\usepackage{natbib}
\usepackage[dvips]{graphicx}

\newtheorem{lem}{Lemma}
\newtheorem{thm}{Theorem}

\def\QED{\mbox{\rule[0pt]{1.5ex}{1.5ex}}}

\def\endproof{\hspace*{\fill}~\QED\par\endtrivlist\unskip}

\def\inner{\mathop{\rm int}\nolimits}

\def\argmax{\mathop{\rm argmax}\nolimits}

\def\real{\mathbb{R}}

\def\Label#1{\label{#1}\ [\ #1\ ]\ }
\def\Label{\label}
\title[Two non-regular extensions of large deviation bound]{Two non-regular extensions of large deviation bound}
\author[Masahito Hayashi]{Masahito Hayashi}
\address{Graduate School of Information Sciences,
Tohoku University, Aoba-ku, Sendai 980-8579, JAPAN}
\email{hayashi@math.is.tohoku.ac.jp}
\begin{document}
\begin{abstract}
We formulate two types of extensions of the large deviation theory initiated by Bahadur in a non-regular setting. 
One can be regarded as a bound of the point estimation, the other can be regarded as the limit of a bound of the interval estimation. Both coincide in the regular case, but do not necessarily coincide in a non-regular case.
Using the limits of relative R\'{e}nyi entropies, 
we derive their upper bounds and give a necessary and sufficient condition 
for the coincidence of the two upper bounds. 
We also discuss the attainability of these two bounds in several non-regular location shift families.
\end{abstract}

\section{Introduction}\Label{s1}
As is known, Fisher information gives the bound of 
the accuracy of estimation.
However, this fact holds only for a regular distribution family.
Hence, if the distribution family does not satisfy the regurality condition,
we have to treat the information quantity as alternative of
Fisher information.
In this paper, we consider this problem only for the location shift family on 
the real line $\real$ such that 
the support depends on the true parameter.
So far, this problem has been mainly studied 
concerning the mean square error.
This paper treats this problem 
with the large deviation criterion,
which was introduced by Bahadur \cite{Ba60,Ba67,Ba71}.
Owing this method, we can discuss the difference between 
the interval estimation and the point estimation.
For this analysis, we introduce two type extensions of
Bahadur's large deviation bound.
One is the quantity $\alpha_1(\theta)$,
which can be regarded as the limit of the accuracy of 
the interval estimation.
The other is the quantity $\alpha_2(\theta)$,
which can be regarded as the limit of the accuracy of the point estimation.
We also show that these two quantities coincide 
in a regular distribution family,
but they do not coincide in a non-regular distribution family.

In order to evaluate 
the actuary of a sequence of estimators $\vec{T}=\{T_n\}$
for a probability distribution family
$\{p_{\theta}|\theta \in \Theta\}$ from the large deviation viewpoint,
Bahadur 
focused on 
the error probability $p^n_\theta\{ | T_n - \theta | \,> \epsilon\}$,
which goes to $0$ exponentially.
That is, the rate of the quantity 
$- \log p^n_\theta\{ | T_n - \theta | \,> \epsilon\}$
is in order $n$.
Hence, he discussed the following exponential decreasing rate 
of the error probability with a fixed error bar \cite{Ba60,Ba67,Ba71}:
\begin{align}
\beta(\vec{T},\theta,\epsilon):=
\liminf \frac{-1}{n} \log p^n_\theta
\{ | T_n - \theta | \,> \epsilon\} \Label{12}.
\end{align}
For this purpose, he treated 
Kullback-Leibler's relative entropy (divergence)
$D(p\|q)
:= \int p(\omega) (\log p(\omega) - \log q(\omega)) \,d \omega$.
Using the relation between its limit and Fisher information,
he characterized 
the slope $\alpha(\vec{T},\theta)
:=
\lim_{\epsilon\to 0}\frac{\beta(\vec{T},\theta,\epsilon)
}{\epsilon^2}$
of the exponential rate by 
the Fisher information.
Since this discussion is a fundamental for this paper,
it is summarized in Section \ref{s2}.

However, this method cannot be applied to 
the family in which the support depends on the true parameter,
because Kullback-Leibler's relative entropy diverges in this case.
On the other hand, it is not impossible to define Fisher information 
even if Kullback-Leibler's relative entropy diverges.
In this paper, we focus on the 
the limit of the relative R\'{e}nyi entropies 
$I^s(p\|q)
:= -\log \int p^s(\omega) q^{1-s}(\omega) \,d \omega ~
(0 \,< s \,< 1)$,
which is deeply treated 
in a non-regular location shift family by Hayashi \cite{Haya-2}.
In Section \ref{s3},
we define two quantities 
$\alpha_1(\theta)$ and $\alpha_2(\theta)$,
and give their upper bounds
$\overline{\alpha}_1(\theta)$ and $\overline{\alpha}_2(\theta)$.
Both of $\alpha_1(\theta)$ and $\alpha_2(\theta)$
are defined as 
the optimal slope of the exponential rate
$\beta(\vec{T},\theta,\epsilon)$
concerning $\epsilon$ at the limit $\epsilon \to 0$,
but their definitions are slightly different.
Also we derive a necessary and sufficient condition 
for the coincidence of these two upper bounds. 
Since these upper bounds are proved by a very general treatment 
in Section \ref{s4},
they are valid for general non-regular families.
In Section \ref{s5}, we focus on 
several estimators useful for
location shift families.
We calculate their exponential rates 
$\beta(\vec{T},\theta,\epsilon)$
and their slopes $\alpha(\vec{T},\theta)$
concerning $\epsilon$ at the limit $\epsilon \to 0$.


\section{Bahadur theory}\Label{s2}
In this section, we begin by summarizing the results reported by Bahadur \cite{Ba60,Ba67,Ba71}, who discussed the decreasing rate of the tail probability in the estimation for a distribution family. Given $n$-i.i.d. data 
$\omega_1, \ldots, \omega_n$, the exponential rate $\beta(\vec{T},\theta,\epsilon)$ of the estimator 
$\vec{T}=\{T_n\}$ is written as 
\begin{align*} \beta(\vec{T},\theta,\epsilon)= \min \{ \beta^+(\vec{T},\theta,\epsilon), \beta^-(\vec{T},\theta,\epsilon)\},
\end{align*}
where the exponential rates of half-side error
probabilities are given by
\begin{align*}
\beta^+(\vec{T},\theta,\epsilon) &:=
\liminf \frac{-1}{n} \log p^n_\theta
\{ T_n \,> \theta + \epsilon\} \\
\beta^-(\vec{T},\theta,\epsilon) &:=
\liminf \frac{-1}{n} \log p^n_\theta
\{ T_n \,< \theta - \epsilon\}.
\end{align*}
When an estimator $\vec{T}=\{T_n\}$ satisfies the weak consistency
\begin{align*}
 p^n_\theta
\{ | T_n - \theta | \,> \epsilon\} \to 0 
\quad \forall \epsilon \,> 0, \quad
\forall \theta \in \Theta,
\end{align*}
using the monotonicity of KL-divergence, we can prove the inequality
\begin{align}
\beta(\vec{T},\theta,\epsilon)
\le \min\{ D(p_{\theta+\epsilon} \|p_\theta),
D(p_{\theta-\epsilon} \|p_\theta)\}. \Label{Baha4}
\end{align}
Note that if, and only if, the family is exponential, there exists an estimator attaining the equality
(\ref{Baha4}) at 
$\forall\theta\in \Theta, \forall\epsilon \,> 0$.
Therefore, for a general family, it is difficult to optimize 
the exponential rate $\beta(\vec{T},\theta,\epsilon)$.

Instead of the exponential rate $\beta(\vec{T},\theta,\epsilon)$,
We usually consider the slope of 
the exponential rate:
\begin{align}
\alpha(\vec{T},\theta)
:= 
\lim _{\epsilon \to + 0}
\frac{1}{\epsilon^2}
\beta(\vec{T},\theta,\epsilon). \Label{2.5}
\end{align}
In this case, when the Fisher information $J_\theta$ satisfies the condition
\begin{align}
J_\theta:= \int l_\theta (\omega)^2 p_\theta(\,d \omega)
=\lim _{\epsilon \to 0}
\frac{2}{\epsilon^2}
D(p_{\theta+\epsilon} \|p_\theta),
\quad l_\theta (\omega):=
\frac{\,d }{\,d \theta}
\log p_{\theta} (\omega),\Label{4-6-10}
\end{align}
the inequality
\begin{align}
\alpha(\vec{T},\theta)
\le
\frac{1}{2} J_\theta
\Label{Baha1}
\end{align}
holds.
Moreover, as was proven by Fu \cite{Fu73}, if the family satisfies the concavity of the logarithmic derivative $l_\theta(\omega)$ for $\theta$ and some other conditions, the MLE $\vec{\theta_{ML}}$ attains the equality of (\ref{Baha1}).
These facts are summarized in the two equations:
\begin{align}
\alpha(\theta) :=
\sup_{\vec{T}:{\rm WC}}
\alpha(\vec{T},\theta) 
&=
\frac{1}{2} J_\theta \Label{Baha2} \\
\lim _{\epsilon \to + 0}
\frac{1}{\epsilon^2}
\sup_{\vec{T}:{\rm WC}}
\beta(\vec{T},\theta,\epsilon)
&=
\frac{1}{2} J_\theta \Label{Baha3}.
\end{align}
As is mentioned later, these equations 
imply an interesting relation 
between the point estimation and
the interval estimation.

\section{Upper bounds}\Label{s3}
In this paper, the relative R\'{e}nyi entropies $I^s(p_\theta\|p_{\theta+\epsilon})$ substitute for the KL divergence. Note that the order of $I^s(p_\theta\|p_{\theta+\epsilon})$ is not necessarily $\epsilon^2$ at the limit $\epsilon \to 0$.
However, its order is independent of the parameter $s$, as is guaranteed by the inequalities
\begin{align*}
2 \min \{ s, 1-s\} I^{\frac{1}{2}}
(p_\theta\|p_{\theta+\epsilon})
\le I^s(p_\theta\|p_{\theta+\epsilon})
\le
2 \max \{ s, 1-s\} I^{\frac{1}{2}}
(p_\theta\|p_{\theta+\epsilon}),
\end{align*}
which are proven in Lemma \ref{L38} of the Appendix. In several cases, 
the order of the exponential rate
$\beta(\vec{T},\theta,\epsilon)$ coincides with the order of the relative R\'{e}nyi entropies $I^s(p_\theta\|p_{\theta+\epsilon})$.
In the following, we use a strictly monotonically decreasing function $g(x)$ such that $I^s(p_\theta\|p_{\theta+\epsilon})\cong O(g(\epsilon)) $ and $g(0)=0$. 

Following equations (\ref{Baha2}) and (\ref{Baha3}), 
we define two extensions of Bahadur's bound (slope) $\alpha(\theta)$ as
\begin{align}
\alpha_1(\theta)&:=
\limsup_{\epsilon \to + 0}
\frac{1}{g(\epsilon)}
\sup_{\vec{T}}
\inf_{\theta-\epsilon \le \theta'\le \theta+\epsilon}
\beta(\vec{T},\theta',\epsilon) \Label{hi1} \\
\alpha_2(\theta)&:=\sup_{\vec{T}}
\alpha_2(\vec{T},\theta),
\Label{hi2}
\end{align}
where
\begin{align}
\alpha_2(\vec{T},\theta):= \liminf_{\epsilon \to + 0}
\frac{1}{g(\epsilon)}
\inf_{\theta-\epsilon \le \theta'\le \theta+\epsilon}
\beta(\vec{T},\theta',\epsilon) .\Label{8-a}
\end{align}
In the interval estimation,
we consider only the error probability concerning the fixed 
interval $[\theta - \epsilon, \theta+\epsilon]$.
That is, the optimization depending on the quantity $\epsilon$
is allowed.
Hence, the performance with enough small $\epsilon$
is characterized by 
$\alpha_1(\theta)$ not by $\alpha_2(\theta)$.
On the other hand,
in the point estimation, 
we have to treat the error probability concerning the 
interval $[\theta - \epsilon, \theta+\epsilon]$ for all $\epsilon$.
The performance of an estimator $\vec{T}$ is 
characterized as the limit $\alpha_2(\vec{T},\theta)$.
Then, the bound of the performance of the point estimation
is given by 
$\alpha_2(\theta):=\sup_{\vec{T}}
\alpha_2(\vec{T},\theta)$.
In the regular family,
equation (\ref{Baha2}) can be regarded as the bound of the point estimation, 
while equation (\ref{Baha3}) can be regarded 
as the limit of the bound of the interval estimation because $\sup_{\vec{T}:WC}\beta(\vec{T},\theta,\epsilon)$ corresponds to the bound of the interval estimation with width $2 \epsilon$ of the confidence interval.
Therefore, we can conclude that 
there is no difference between the point estimation and
the limit of the interval estimation in the estimation in the regular 
family.
In the following, we consider 
whether there exists a difference between them.

Note that we take infimum
$\inf_{\theta-\epsilon \le \theta'\le \theta+\epsilon}$
into account in (\ref{8-a}), unlike (\ref{2.5}).
As was pointed out by Ibragimov and Has'minskii
\cite{IH},
when KL-divergence is infinite, there exists a consistent super efficient estimator $\vec{T}$ such that $\beta(\vec{T},\theta,\epsilon)$ and $\lim _{\epsilon \to + 0}
\frac{1}{g(\epsilon)} \beta(\vec{T},\theta,\epsilon)$ is infinite at one point $\theta$.
Therefore, we need to take the infimum 
$\inf_{\theta-\epsilon \le \theta'\le \theta+\epsilon}$
into account. In this situation, we do not need to 
limit estimators to weakly consistent ones. As is proven in the next section, we can obtain the following theorems. 
\begin{thm}\Label{thm1}
When the convergence $\lim_{\epsilon \to 0 } \frac{I^s(p_{\theta-\epsilon/2}\| p_{\theta+\epsilon/2})} {g(\epsilon)}$ 
is uniform for $0\,< s \,< 1$, the inequality
\begin{align}
\alpha_1(\theta) \le 
\overline{\alpha}_1(\theta):= 2^{\kappa}
\sup_{0 \,< s \,< 1}I^s_{g,\theta} \Label{7-a}
\end{align}
holds, where $\kappa$ and 
$I^s_{g,\theta}$ are defined by
\begin{align}
I^s_{g,\theta} &:= \lim_{\epsilon \to +0 } 
\frac{I^s(p_{\theta-\epsilon/2}\| 
p_{\theta+\epsilon/2})}
{g(\epsilon)} \quad 1 \ge s \ge 0 \Label{uni} \\
x^\kappa &= \lim_{\epsilon \to +0} \frac{g(x\epsilon)}
{g(\epsilon)} . \Label{4-6}
\end{align}
Lemma \ref{ap1} proven in Appendix \ref{A}, 
guarantees the existence of such a real number $\kappa$.
\end{thm}
Note that the function $s \to I^s_g$ is concave and continuous because the function $s \to I^s(p_{\theta-\frac{1}{2} \epsilon}\|p_{\theta+\frac{1}{2}\epsilon})$ is concave and continuous.
Therefore, when $I_{g,\theta}^s=I_{g,\theta}^{1-s}$, we have
\begin{align}
\overline{\alpha}_1(\theta)= 2^\kappa
I^{\frac{1}{2}}_{g,\theta}.
\end{align}
\begin{thm}\Label{thm2}
If the convergence $\lim_{\epsilon \to 0 } \frac{I^s(p_{\theta-\epsilon/2}\| p_{\theta+\epsilon/2})} {g(\epsilon)}$ is uniform for $ s \in (0,1)$ and $\theta \in K$ for any compact set $K \subset \real$, the inequality
\begin{align}
\alpha_2(\theta)
\le
\overline{\alpha}_2(\theta) 
:=
\left\{
\begin{array}{cc}
\sup_{0\,< s \,< 1}
\frac{I^s_{g,\theta}}{s(1-s)}
\left( s^{\frac{1}{\kappa-1}} + (1-s)^{\frac{1}{\kappa-1}} 
\right)^{\kappa-1}& \hbox{ if }
\kappa \,< 1 \\
2 I^{\frac{1}{2}}_{g,\theta} & \hbox{ if }
\kappa = 1 \\
\inf_{0 \,< s \,< 1}
\frac{I^s_{g,\theta}}{s(1-s)}
\left( s^{\frac{1}{\kappa-1}} + (1-s)^{\frac{1}{\kappa-1}} 
\right)^{\kappa-1}& \hbox{ if }
\kappa \,> 1 
\end{array}
\right. 
\Label{9-a}
\end{align}
holds.
\end{thm}
In our proofs of these theorems, Chernoff's formula and Hoeffding's formula in simple hypothesis testing play important roles.

As was proven by Akahira \cite{Ak}, under some regularity conditions for a distribution family, the equation
\begin{align}
\lim_{\epsilon \to +0}
\frac{I^s(p_\theta\|p_{\theta+\epsilon})}
{\epsilon^2}
= \frac{J_\theta s(1-s)}{2}, \Label{38}
\end{align}
holds.
When we choose $g(x)=x^2$, we have $\kappa=2, I^s_{g,\theta} = \frac{1}{2}J_{\theta}s(1-s)$, and the relations
\begin{align*}
\overline{\alpha}_1(\theta) &= 
4 \max_{0\le s \le 1}
I^s_{g,\theta}
=\frac{1}{2} J_\theta \\
\overline{\alpha}_2(\theta) &= 
\min_{0\le s \le 1}
\frac{I^s_{g,\theta}}{s(1-s)}
=\frac{1}{2} J_\theta 
\end{align*}
hold.
In particular, if the distribution family satisfies the concavity of the logarithmic derivative 
$l_\theta(\omega)$ for $\theta$ and some other conditions, the bound $\frac{1}{2} J_\theta $ is attained by the MLE. Thus, 
the relations $\alpha_1(\theta) = \alpha_2(\theta) = \overline{\alpha}_1(\theta) = \overline{\alpha}_2(\theta) 
= \frac{1}{2} J_\theta$ hold.

As a relation between two bounds $\overline{\alpha}_1
(\theta) $ and $\overline{\alpha}_2(\theta)$, we can 
prove the following theorem.
\begin{thm}\Label{thm3}
The inequality
\begin{align}
\overline{\alpha}_1(\theta) 
\ge 
\overline{\alpha}_2(\theta)  \Label{14}
\end{align}
holds, and (\ref{14}) holds as an equality 
if and only if the equations 
\begin{align}
\overline{\alpha}_1(\theta) 
&=
2^\kappa I_{g,\theta}^{\frac{1}{2}} \Label{16.3} \\
2^\kappa I_{g,\theta}^{\frac{1}{2}}
&=
\overline{\alpha}_2(\theta) \Label{15}
\end{align}
hold.
When $\kappa \le 1$, 
(\ref{14}) holds as an equality 
if and only if equation (\ref{16.3}) holds. 
\end{thm}
When $I^s_{g,\theta}$ is differentiable, condition (\ref{16.3}) is equivalent to $ \left. \frac{\,d }{\,d s} {I}_{g,\theta}^{s}
\right|_{s=\frac{1}{2}} =0$.

\section{Proofs of upper bounds}\Label{s4}
In our proofs of Theorems \ref{thm1} and \ref{thm2}, 
Chernoff's formula and Hoeffding's formula in simple 
hypothesis testing are essential, and are summarized 
as follows. 
Let the probability $p$ on $\Omega$ be the null 
hypothesis and $q$ be the alternative hypothesis. 
When we discuss 
a hypothesis testing problem concerning $n$-i.i.d. 
data, we call a sequence $\vec{A}=\{A_n\}$ a test, 
where $A_n$ is an acceptance region, which is a subset 
of $\Omega^n$. The first error probability $e_1(A_n)$ 
and the second error probability $e_2(A_n)$ are defined as
\begin{align*}
e_1(A_n):= 1- p^n(A_n) , \quad
e_2(A_n):= q^n(A_n),
\end{align*}
and their exponents are given by
\begin{align*}
e^*_1(\vec{A}):= \liminf_{n\to \infty}
\frac{-1}{n} \log e_1(A_n) \\
e^*_2(\vec{A}):= \liminf_{n\to \infty}
\frac{-1}{n} \log e_2(A_n) .
\end{align*}
Chernoff \cite{Cher} evaluated the exponent of the sum of the two errors as 
\begin{align}
\sup_{\vec{A}}
\lim_{n \to \infty} \frac{-1}{n}
\log
(e_1(\vec{A}_n)+e_2(\vec{A}_n))
=
\sup_{\vec{A}}
\min\{ e^*_1(\vec{A}),
e^*_2(\vec{A})\}
=\sup_{0\,< s \,< 1} I^s(p\|q),  \Label{11}
\end{align}
which is essential for our proof of Theorem \ref{thm1}.
This bound is achieved by both of the likelihood tests $\{ \vec{\omega}_n \in \Omega^n| p^n(\vec{\omega}_n) \ge q^n(\vec{\omega}_n)\}$ and $\{ \vec{\omega}_n \in \Omega^n| p^n(\vec{\omega}_n) \,> q^n(\vec{\omega}_n)\}$.

Hoeffding proved another formula for simple hypothesis testing \cite{Hoef}:
\begin{align}
\sup_{\vec{A}:e^*_1(\vec{A})\ge r }
e^*_2(\vec{A}) = \sup_{0 \,< s \,< 1}
\frac{-sr + I^s(p\|q)}{1-s}.
 \Label{Hoe}
\end{align}
This formula is essential for our proof of Theorem \ref{thm2}.

\noindent\hspace{2em}{\it Proof of Theorem \ref{thm1}: }
Applying equation (\ref{11}) to the two hypotheses $p_{\theta-\epsilon}$ and $p_{\theta+\epsilon}$, we obtain 
\begin{align*}
\inf_{\theta-\epsilon\le \theta'
\le \theta+\epsilon}
\beta(\vec{T},\theta',\epsilon)\le
\min\{ \beta^+(\vec{T},\theta-\epsilon,\epsilon),
\beta^-(\vec{T},\theta+\epsilon,\epsilon) \}
\le
\sup_{0 \,< s \,< 1} I^s(p_{\theta-\epsilon} \|
p_{\theta+\epsilon}).
\end{align*}
Taking the limit $\epsilon \to 0$, we have
\begin{align}
\lim _{\epsilon \to + 0}
\frac{1}{g(\epsilon)}
\sup_{\vec{T}:WC} 
\inf_{\theta-\epsilon\le \theta'
\le \theta+\epsilon}
\beta(\vec{T},\theta',\epsilon)
\le
2^{\kappa}
\sup_{0 \,< s \,< 1}I^s_{g,\theta}. \Label{hise1}
\end{align}
\endproof
In the regular case, this method was used by Sievers \cite{Sie}.

\noindent\hspace{2em}{\it Proof of Theorem \ref{thm2}:} 
Hoeffding's formula (\ref{Hoe}) yields
\begin{align*}
& \inf_{\theta-(1-\eta)\epsilon \le 
 \theta'\le \theta+ (1-\eta)\epsilon}
\beta(\vec{T},\theta',(1-\eta)\epsilon) \\
\le &
\beta^-(\vec{T},\theta+ \frac{1}{2}\epsilon,(1-\eta)\epsilon) \\
\le &
\sup_{0\,< s \,< 1}
\frac{-s 
\beta^+(\vec{T},\theta- \frac{1}{2}\epsilon,\eta\epsilon)
+ I^s(p_{\theta- \frac{1}{2}\epsilon}\|
p_{\theta+ \frac{1}{2}\epsilon})}
{1-s} \\
\le &
\sup_{0\,< s \,< 1}
\frac{-s 
\inf_{\theta-\eta\epsilon \le \theta'\le \theta+ \eta\epsilon}
\beta(\vec{T},\theta',\eta\epsilon)
+ I^s(p_{\theta- \frac{1}{2}\epsilon}\|
p_{\theta+ \frac{1}{2}\epsilon})}{1-s}.
\end{align*}
The uniformity of (\ref{uni}) guarantees that
\begin{align}
\alpha_2(\vec{T},\theta)
(1-\eta)^\kappa
&\le
\sup_{0\,< s \,< 1}
\frac{-s \alpha_2(\vec{T},\theta) \eta^\kappa
+ I^s_{g,\theta}}{1-s}.
\Label{1-a}
\end{align} From (\ref{1-a}), we have
\begin{align}
\alpha_2(\vec{T},\theta)
\le 
\sup\left\{
\alpha\left|
\alpha ( \eta^\kappa, (1-\eta)^\kappa)
\in
\left\{
(x,y)\left|
y \le \max_{0 \le s\le 1}
\frac{-sx +I^s_g}
{1-s}, x, y \ge 0
\right.
\right\}, \quad
0 \le \forall \eta \le 1
\right.
\right\}. \Label{cut3}
\end{align}
We define the set ${\cal C}_1$ and $\alpha_0$ as
\begin{align}
{\cal C}_1&:=
\left\{
(x,y)\left| y \ge \sup_{0\,< t \,< 1}
\frac{-tx +I^t_{g,\theta}}{1-t}
\right.\right\} \\
\alpha_0&:= \sup\{\alpha|
(\alpha \eta^{\kappa},
\alpha(1-\eta)^{\kappa})
\notin \inner({\cal C}_1),
\quad 0\le \forall \eta \le 1\}. \Label{28-1}
\end{align}
Note that the function $s \mapsto I^s_{g,\theta}$ is concave. We define the convex function $g(x)$ and another set ${\cal C}_2$ as
\begin{align}
g(x)&:=\left\{
\begin{array}{cc}
\alpha_0 \left( 1- \left(\frac{x}{\alpha_0}\right)
^{\frac{1}{\kappa} }\right)^{\kappa}&
0\le x \le \alpha_0 \\
0 & x \,> \alpha_0 
\end{array}
\right. \Label{28-2}\\
{\cal C}_2 &:= \{ (x,y) | y \ge g(x)\}.\Label{28-3}
\end{align}
Since inequality (\ref{cut3}) guarantees 
\begin{align*}
(\alpha_2(\vec{T},\theta)\eta^{\kappa},
\alpha_2(\vec{T},\theta)(1-\eta)^{\kappa})
\notin \inner({\cal C}_1), 
\quad 0\le \forall \eta \le 1,
\end{align*}
the inequality 
\begin{align}
\alpha_2(\vec{T},\theta) \le \alpha_0
\end{align}
holds. Relations (\ref{28-1}), (\ref{28-2}), and (\ref{28-3})
guarantee the 
relation
\begin{align}
{\cal C}_1 \subset {\cal C}_2 \Label{28-4}.
\end{align}
Applying Lemma \ref{concave}, we have
\begin{align}
& \inf_{(x,y)\in {\cal C}_1}
(sx +(1-s)y)
=
\inf_{x \ge 0}
\left(sx +(1-s)\sup_{0 \,< t \,< 1}
\frac{-tx +I^t_{g,\theta}}{1-t}\right) \nonumber\\
=&
\inf_{x \,>0}\sup_{0 \,< t \,< 1}
\left( \frac{(s-t)x +(1-s)I^t_{g,\theta}}{1-t}\right)
= I^s_{g,\theta}. \Label{4-19}
\end{align}

In the following, we divide our situation into three cases $\kappa \,> 1, \kappa =1 , 1 \,> \kappa \,> 0$. When $\kappa \,> 1$, relation (\ref{28-4}) guarantees that 
\begin{align*}
&\alpha_0 s(1-s)\left( s^{\frac{1}{\kappa-1}}
+ (1-s)^{\frac{1}{\kappa-1}}\right)
^{1-\kappa}
=
\alpha_0 \min_{0\,< \eta \,< 1}
\left(s\eta^{\kappa}+(1-s)(1-\eta)^{\kappa}\right)\\
=& \inf_{(x,y) \in {\cal C}_2} \left(sx + (1-s)y \right)
\le  \inf_{(x,y) \in {\cal C}_1} (sx + (1-s)y )
= I^s_{g,\theta}.
\end{align*}
Therefore, 
\begin{align*}
\alpha_2(\vec{T},\theta)
\le
\alpha_0
\le
\frac{ \left( s^{\frac{1}{\kappa-1}}
+ (1-s)^{\frac{1}{\kappa-1}}\right)
^{\kappa-1}}{s(1-s)}
I^s_{g,\theta},
\end{align*}
which implies (\ref{9-a}).

When $\kappa =1$, similarly, we can easily prove 
\begin{align*}
\frac{1}{2}
\alpha_2(\vec{T},\theta)
=\alpha_2(\vec{T},\theta)
\min_{0\,< \eta\,<1} \left(\frac{1}{2}\eta
+\frac{1}{2}(1-\eta)\right)
\le \inf_{(x,y)\in {\cal C}_1} 
\left(\frac{1}{2} x + \frac{1}{2} y \right)
= I^{\frac{1}{2}}_{g,\theta}.
\end{align*}

Finally, we consider the case where $\kappa \,< 1$. Since function $g$ is concave on $(0,\alpha_0)$, there exists $\eta_0 \in (0,1)$ such that
\begin{align}
\alpha_0 (\eta_0^{\kappa} ,(1-\eta_0)^{\kappa})
\in {\cal C}_1 \cap \overline{{\cal C}_2^c}.
\end{align}
Since ${\cal C}_1$ and ${\cal C}_2^c$ are convex, there exists a real number $s_0 \in (0,1)$ such that 
\begin{align*}
\sup_{(x,y) \in  {\cal C}_2^c}
\left(s_0 x + (1-s_0)y \right)
=
\alpha_0 ( s_0 \eta_0^{\kappa} 
+ (1-s_0)(1-\eta_0)^{\kappa})
=\inf_{(x,y) \in  {\cal C}_1}
\left(s_0 x + (1-s_0)y\right).
\end{align*}
In general, for any $s \in (0,1)$, using (\ref{4-19}), we obtain 
\begin{align*}
& \alpha_0 s(1-s) \left(s^{\frac{1}{\kappa-1}}
+(1-s)^{\frac{1}{\kappa-1}}\right)
^{1-\kappa} 
= \alpha_0 \sup_{0\,< \eta \,< 1}
\left(s \eta^{\kappa} +(1-s) (1-\eta)^{\kappa} \right)\\
=&\sup_{(x,y) \in  {\cal C}_2^c}
\left(s x + (1-s)y \right)
\ge  \alpha_0 \left(
s \eta_0^{\kappa} +(1-s) (1-\eta_0)^{\kappa}\right)
\ge 
\inf_{(x,y) \in  {\cal C}_1}
\left(s x + (1-s)y\right) =
I^{s}_{g,\theta},
\end{align*}
which lead to (\ref{9-a}).
\endproof

\noindent\hspace{2em}{\it Proof of Theorem \ref{thm3}:}
It is trivial in the case of $\kappa=1$.
When $\kappa \,>1$, it is also trivial because 
\begin{align*}
\inf_{0 \,< s \,< 1}
\frac{I^s_{g,\theta}}{s(1-s)}
\left( s^{\frac{1}{\kappa-1}} 
+ (1-s)^{\frac{1}{\kappa-1}} 
\right)^{\kappa-1}
\le
2^\kappa I_{g,\theta}^{\frac{1}{2}}
\le
2^\kappa \sup_{0\le s \le 1}
I_{g,\theta}^{s}.
\end{align*}
Next, we consider the case $\kappa \,< 1$. The inequality
\begin{align}
\sup_{0\,< s \,< 1}
\frac{I^s_{g,\theta}}{s(1-s)}
\left( s^{\frac{1}{\kappa-1}} 
+ (1-s)^{\frac{1}{\kappa-1}} 
\right)^{\kappa-1}
\le
2^\kappa \sup_{0\le s \le 1}
I_{g,\theta}^{s} \Label{19}
\end{align}
follows from the two inequalities
\begin{align}
I^s_{g,\theta}
&\le \sup_{0 \,< s \,< 1}
I_{g,\theta}^{s} \Label{17}\\
\frac{1}{s(1-s)}
\left( s^{\frac{1}{\kappa-1}} 
+ (1-s)^{\frac{1}{\kappa-1}} 
\right)^{\kappa-1}
&\le
2^\kappa \Label{18.9}.
\end{align}
We thus obtain (\ref{14}). 
In the following, we prove that the equality of (\ref{14}) implies (\ref{16.3}) and (\ref{16.3}) implies (\ref{15}) and the equality of (\ref{14}) in the case where $\kappa \,< 1$. If we assume that the equality of (\ref{19}) holds, the equalities of (\ref{17}) and (\ref{18.9}) hold at the same $s$. The equality of (\ref{18.9}) holds 
if and only if $s=\frac{1}{2}$. Therefore, 
\begin{align*}
I^{\frac{1}{2}}_{g,\theta}
= \sup_{0 \,< s \,< 1}
I_{g,\theta}^{s},
\end{align*}
which is equivalent to (\ref{16.3}). If we assume that (\ref{16.3}) holds, inequality (\ref{18.9}) guarantees that
\begin{align}
\overline{\alpha}_2(\theta) 
\le 2^{\kappa} I^{\frac{1}{2}}_{g,\theta}.
\end{align}
Substituting $\frac{1}{2}$ into $s$ at the 
left hand side (LHS) in the definition of $\overline{\alpha}_2(\theta)$, we obtain
\begin{align*}
\overline{\alpha}_2(\theta) \ge
2^{\kappa} I^{\frac{1}{2}}_{g,\theta}.
\end{align*}
Thus, equation (\ref{15}) holds. Combining (\ref{16.3}) and (\ref{15}), we obtain the equality of (\ref{14}).
\endproof
\section{Exponential rates of useful estimators and their slopes}
\Label{s5}
In the following, we discuss the exponential rate 
$\beta(\vec{T}, \theta,\epsilon)$ of a useful estimator $\vec{T}$ 
for a location shift family $\{ f(x - \theta)| \theta \in \real\}$, 
where $f$ is a probability density function (pdf) on $\real$.
In particular, we focus on the case where the support 
of $f$ is $(a,b)$.
Further, we assume that the pdf $f$ is $C^1$ continuous and satisfies that
\begin{align*}
f(x) \cong A_1 (x-a)^{\kappa_1-1}, & \quad x \to a +0\\
f(x) \cong A_2 (b- x)^{\kappa_2-1},& \quad x \to b -0,
\end{align*}
where $\kappa_1,\kappa_2 \,> 0$, as for the beta distributions
$\frac{x^{\alpha-1}(1-x)^{\beta-1}}{B(\alpha,\beta)}$.

When its support is a half line $(0,\infty)$ as for the gamma distribution and Weibull distribution, our situation results in the above case where $A_2=0$ if $f$ is $C^3$ continuous and $\lim_{x \to \infty}
|\frac{\,d }{\,d x} \log f(x)|\,< \infty$.
Also, when $\kappa_1 \,> \kappa_2$, our situation results in the above case where $A_2=0$.

In the first step, we will treat estimators useful for point estimation.
After this discussion, we will discuss estimators useful for
interval estimation.
\subsection{Exponential rates of estimators useful for point estimation}
When the support of $f$ is $(a,b)$, 
the two estimators $\overline{\theta}_n:= \max\{ x_1, \ldots, x_n \} - b$ and $\underline{\theta}_n:= \min\{ x_1, \ldots, x_n \} - a$
are useful for point estimation.
These performances are characterized as follows.
\begin{lem}\Label{l34.3}
The estimators $\vec{\underline{\theta}}$ and $\vec{\overline{\theta}}$ satisfy 
\begin{align}
\beta^+(\vec{\underline{\theta}},\theta,\epsilon)
&=
- \log \left(\int_{a}^{b-\epsilon}f(x) \,d x\right)
,\quad
\beta^-(\vec{\underline{\theta}},\theta,\epsilon)
=\infty \Label{3-8-14} \\
\beta^+(\vec{\overline{\theta}},\theta,\epsilon)
&=\infty ,\quad
\beta^-(\vec{\overline{\theta}},\theta,\epsilon)
=
- \log \left(\int_{a+\epsilon}^{b}f(x) \,d x \right).
\Label{3-8-15}
\end{align}
\end{lem}
\begin{proof}
Since 
\begin{align*}
f_\theta^n\{\underline{\theta}_n>\theta+\epsilon\}
= \left(\int_{a}^{b-\epsilon}f(x) \,d x\right)^n ,\quad 
f_\theta^n\{\underline{\theta}_n <\theta\}=0\\
f_\theta^n\{\overline{\theta}_n<\theta-\epsilon\}
= \left(\int_{a+\epsilon}^{b}f(x) \,d x \right)^n,\quad
f_\theta^n\{\overline{\theta}_n>\theta\}=0,
\end{align*}
we obtain (\ref{3-8-14}) and (\ref{3-8-15}).
\end{proof}
In order to strike a balance between 
two exponential rates 
$\beta^+(\vec{T},\theta,\epsilon)$ and 
$\beta^-(\vec{T},\theta,\epsilon)$,
we define the convex combination (CC) estimator 
$\vec{\theta}_{CC,\lambda}:= 
\{ \theta_{CC,\lambda,n}:= 
\lambda \underline{\theta}_n + (1-\lambda) \overline{\theta}_n\}$ 
with the ratio $\lambda : 1- \lambda$
of the two estimators $\vec{\underline{\theta}}$ and $\vec{\overline{\theta}}$, where $0 \,< \lambda \,< 1$. 
These are characterized as follows.
\begin{lem}\Label{l34.3-2}
The convex combination (CC) estimator 
$\vec{\theta}_{CC,\lambda}$ satisfies that 
\begin{align}
\beta^+(\vec{\theta}_{CC,\lambda},\theta,\epsilon)
&=
- \log \left(\int_{a}^{b-\frac{\epsilon}{1-\lambda}}f(x) \,d x\right)
\Label{34.1}\\
\beta^-(\vec{\theta}_{CC,\lambda},\theta,\epsilon)
&=
- \log \left(\int_{a+\frac{\epsilon}{\lambda}}^{b}f(x) \,d x \right)
\Label{34.2}.
\end{align}
\end{lem}
\begin{proof}
Define $\overline{\omega}_n:=\max \{ \omega_1 , \ldots , \omega_n \}, \underline{\omega}_n:=\min \{ \omega_1 , \ldots , \omega_n \}$. 
Since the estimators $\vec{\underline{\theta}}, \vec{\overline{\theta}}$, and 
$\vec{\theta}_{CC,\lambda}$ 
are covariant for location shift, we may discuss only the case that $\theta=0$. From the relation $\underline{\theta}_n \,> \theta \,> \theta - \epsilon$, we obtain the second equation of (\ref{3-8-14}). Its joint probability density function $f_n( \overline{\omega}_n,\underline{\omega}_n)$ is given by
\begin{align}
f_n( \overline{\omega}_n,\underline{\omega}_n)
:=
\begin{cases}
n(n-1)\left( \int_{\underline{\omega}_n}^{\overline{\omega}_n}
f(x) \,d x \right)^{n-2}f(\underline{\omega}_n)f(\overline{\omega}_n)
& \overline{\omega}_n \ge \underline{\omega}_n \\
0& \overline{\omega}_n \,< \underline{\omega}_n .
\end{cases}\Label{4-5-1}
\end{align}
By defining
\begin{align*}
g(\underline{\omega}_n,\overline{\omega}_n):=
\begin{cases}
\int_{\underline{\omega}_n}^{\overline{\omega}_n}
f(x) \,d x & \overline{\omega}_n \ge \underline{\omega}_n \\
0& \overline{\omega}_n \,< \underline{\omega}_n ,
\end{cases}
\end{align*}
The equation (\ref{4-5-1}) yields
\begin{align}
p^n_{\theta}(\theta_{CC,\lambda,n} \le \theta- \epsilon)
= \int_{ \theta_{CC,\lambda,n}
(\underline{\omega}_n, \overline{\omega}_n) \le \epsilon}
n(n-1) g( \underline{\omega}_n,\overline{\omega}_n)^n
f(\underline{\omega}_n)f(\overline{\omega}_n)
\,d \underline{\omega}_n
\,d \overline{\omega}_n \Label{34.6}.
\end{align} From the continuity of 
$f(\underline{\omega}_n), f(\overline{\omega}_n)$ and
$g( \underline{\omega}_n,\overline{\omega}_n)$,
the equations 
\begin{align*}
\lim_{n \to \infty}
\frac{1}{n}\log p^n_{\theta}(\theta_{CC,\lambda,n}
\le \theta- \epsilon)
&= 
\lim_{n \to \infty}\frac{1}{n}
\sup_{\theta_{CC,\lambda,n}(\underline{\omega}_n, \overline{\omega}_n)  \le \epsilon}
\log g(\underline{\omega}_n,\overline{\omega}_n) \\
&= - \log \left(\int_{a}^{b-\frac{\epsilon}{1-\lambda}}f(x) \,d x\right)
\end{align*}
hold.
This implies the first equation of (\ref{3-8-14}) and (\ref{34.1}).
In addition, we can similarly show the same fact for (\ref{34.2}).
\end{proof}
Next, we focus on the maximum likelihood estimator $\vec{\theta}_{ML}
:=\{\theta_{ML,n}\}$.

\begin{lem}\Label{l10}
When the function $x \mapsto \log f(x)$ is concave, MLE $\vec{\theta}_{ML}$ 
satisfies that 
\begin{align}
\beta^+(\vec{\theta}_{ML}, \theta,\epsilon) 
&= \sup_{t \ge 0}
- \log \int_{a+\epsilon}^b \exp \left( -t 
\frac{f'(x-\epsilon)}{f(x-\epsilon)} \right)
f(x) \,d x \Label{33} \\
&= \sup_{t \ge 0}
- \log \int_{a}^{b-\epsilon} \exp \left( -t 
\frac{f'(x)}{f(x)} \right)
f(x+\epsilon) \,d x \Label{34}\\
\beta^-(\vec{\theta}_{ML}, \theta,\epsilon) 
&= \sup_{t \ge 0}
- \log \int_{a}^{b-\epsilon} \exp \left( t 
\frac{f'(x+\epsilon)}{f(x+\epsilon)} \right)
f(x) \,d x \Label{35}\\
&= \sup_{t \ge 0}
- \log \int_{a+\epsilon}^b \exp \left( t 
\frac{f'(x)}{f(x)} \right)
f(x-\epsilon) \,d x \Label{36}.
\end{align}
Also, we have the following evaluations in the same assumption.
\begin{align}
\beta^+(\vec{\theta}_{ML}, \theta,\epsilon) 
&\ge \sup_{0 \le s \le 1}I^s(f_{\theta}\|f_{\theta+\epsilon})
\Label{4-5-3} \\
\beta^-(\vec{\theta}_{ML}, \theta,\epsilon) 
&\ge \sup_{0 \le s \le 1}I^s(f_{\theta-\epsilon}\|f_{\theta})\Label{4-5-4}.
\end{align}
\end{lem}
This lemma is essentially proved as a one step of the proof of the main theorem in 
the paper \cite{Fu73}. 
However, he proved the main theorem with a more general assumption and a different notation.
Hence, it is not easy to find the relationship between the notation of the present paper and that of Fu's paper.
Further, since his proof of this part is too short,
it seems that a non-expert of large deviation theory cannot follow his proof.
Therefore, for the reader's convenience, a proof of this lemma is given as follows.

\begin{proof}
Equations (\ref{34}) and (\ref{36}) are trivial. We prove (\ref{33}). From the assumption that for any $\vec{x}_n:= (x_1 , \ldots, x_n)$, the function $\theta \mapsto \sum_{i =1}^n \log f(x_i - \theta)$ is concave 
on $(\overline{\theta}(\vec{x}_n), \underline{\theta}(\vec{x}_n))$, the function $\theta \mapsto \sum_{i=1}^n l_\theta(x_i)$ is monotonically decreasing, 
where $l_\theta(x):=-\frac{f'(x-\theta)}{f(x-\theta)}$ on $(\overline{\theta}(\vec{x}_n), \underline{\theta}(\vec{x}_n))$. 
When $\vec{x}_n$ belongs to the support of $f_{\theta'}$, the condition $\theta_{ML,n}(\vec{x}_n) \ge \theta'$ is equivalent to the condition 
\begin{align*}
\frac{1}{n} \sum_{i=1}^n l_{\theta'}(x_i) \ge 0.
\end{align*}
Denoting the conditional probability $f\{ A| x \in B\}$ under the condition $x \in B$, we can evaluate 
\begin{align}
&\lim - \frac{1}{n} \log f_\theta^n 
\{ \theta_{ML,n} \ge \theta+\epsilon\}\nonumber\\
=&
\lim -\frac{1}{n}\log f_{\theta}^n \{
\theta_{ML,n} \ge \theta+\epsilon |
\vec{x}_n \in (a+\epsilon,b+\epsilon)^n\}
- \frac{1}{n}
\log f^n_{\theta}(a+\epsilon,b+\epsilon)^n 
\nonumber\\
=& \lim -\frac{1}{n}\log f_{\theta,\epsilon}^n 
\left \{
\frac{1}{n} \sum_{i=1}^n l_{\theta+\epsilon}(x_i) 
\ge 0 \right\}
-\log \int_{a+\epsilon}^b f(x) \,d x , \Label{31}
\end{align}
where the probability density function $f_{\theta,\epsilon}$
is defined on the support $(a+\epsilon,b)$ by
\begin{align*}
f_{\theta,\epsilon}(x):=
\frac{f(x)}{\int_{a+\epsilon}^b f(x) \,d x}.
\end{align*}
Chernoff's theorem (Theorem 3.1 in Bahadur \cite{Ba71}) guarantees that
\begin{align} 
& \lim -\frac{1}{n}\log f_{\theta,\epsilon}^n 
\{
\frac{1}{n} \sum_{i=1}^n l_{\theta+\epsilon}(x_i) 
\ge 0 \}\nonumber\\
=&
\sup_{t \ge 0}
- \log \int_{a+\epsilon}^b \exp \left(
t l_{\theta+\epsilon}(x) \right)
f_{\theta,\epsilon}(x) \,d x \nonumber\\
=&
\sup_{t \ge 0}
- \log \int_{a+\epsilon}^b \exp \left(
-t \frac{f'(x-\theta)}{f(x-\theta)}\right)
f(x) \,d x
+\log \int_{a+\epsilon}^b f(x) \,d x .\Label{32}
\end{align}
Combining (\ref{31}) and (\ref{32}), we obtain (\ref{33}). Similarly, 
we can prove (\ref{35}).

Next, in order to show (\ref{4-5-3}),
we choose a real number $\delta > 0$.
When 
\begin{align}
f_{\theta'-\delta}^n (\vec{x}_n) \le f_{\theta'}^n (\vec{x}_n) ,
\end{align}
the condition 
\begin{align}
\theta_{ML,n}(\vec{x}_n) \ge \theta'
\end{align}
holds.
Hence, substituting $\theta-\epsilon $ and $\epsilon$ into 
$\theta'$ and $\delta$, respectively,
we have
\begin{align}
f_{\theta}^n
\{\theta_{ML,n}(\vec{x}_n) \ge \theta+\epsilon \}
\le 
f_{\theta}^n
\{
f_{\theta}^n (\vec{x}_n) \le f_{\theta+\epsilon} (\vec{x}_n) 
\}.
\end{align}
As is mentioned at (\ref{11}) in Section \ref{s4}, 
the exponential rate of the RHS is 
equal to 
$\sup_{0 \le s \le 1}I^s(f_{\theta}\|f_{\theta+\epsilon})$.
Thus, we obtain (\ref{4-5-3}).
Similarly, we can prove (\ref{4-5-4}).
\end{proof}
\begin{lem}\Label{4-5-7}
When $f(x)$ is monotonically decreasing, the MLE $\theta_{ML,n}$ equals the estimator $\underline{\theta}_n$.
\end{lem}
\begin{proof}
For any data $\vec{x}_n:= (x_1, \ldots, x_n)$, if $\theta \,< \underline{\theta}_n(\vec{x}_n))$, $f(x_1- \theta) \ldots f(x_n- \theta)=0$. Conversely, if $\theta \,> \underline{\theta}_n(\vec{x}_n))$, the obtained $f(x_i- \theta) \le f(x_i - \underline{\theta}_n(\vec{x}_n))$. Thus, $\underline{\theta}_n$ is the MLE. 
\end{proof}

\subsection{Slopes of exponential rates of estimators useful for point estimation}
We discuss the slopes of exponential rates
of estimators discussed in the above.
First, we focus on the estimators 
$\vec{\underline{\theta}}$ and $\vec{\overline{\theta}}$.
From Lemma \ref{l34.3},
the estimators 
$\vec{\underline{\theta}}$ and $\vec{\overline{\theta}}$
satisfy that
\begin{align}
\beta(\vec{\underline{\theta}},\theta,\epsilon)
\cong\frac{A_1}{\kappa_1}\epsilon^{\kappa_1}+ o(\epsilon^{\kappa_1}),\quad
\beta(\vec{\overline{\theta}},\theta,\epsilon)
\cong\frac{A_2}{\kappa_2}\epsilon^{\kappa_2}+o(\epsilon^{\kappa_2}).
\Label{4-5-6}
\end{align}

Next, we proceed to the convex combination estimator
$\vec{\theta}_{CC,n}$.
When $\kappa_1=\kappa_2=\kappa$,
Lemma \ref{l34.3-2} yields 
the equations 
\begin{align*}
\beta^+(\vec{\theta}_{CC,\lambda}, \theta, \epsilon )
\cong
\frac{A_1}{\kappa \lambda} \epsilon^{\kappa}+ o(\epsilon^{\kappa})
, \quad
\beta^-(\vec{\theta}_{CC,\lambda}, \theta, \epsilon )
\cong
\frac{A_2}{\kappa(1-\lambda)} \epsilon^{\kappa}
+o(\epsilon^{\kappa}).
\end{align*}
Since $\lambda_0 := 
\frac{A_1^{\frac{1}{\kappa}}}
{A_1^{\frac{1}{\kappa}}+A_2^{\frac{1}{\kappa}}}
= \argmax_{0 \le \lambda \le 1}
\min\{ \frac{A_1}{\kappa \lambda},\frac{A_2}{\kappa(1-\lambda)}\}$,
the relations
\begin{align}
\beta^+(\vec{\theta}_{CC,\lambda_0}, \theta, \epsilon )
\cong
\beta^-(\vec{\theta}_{CC,\lambda_0}, \theta, \epsilon )
\cong
\frac{(A_1^{\frac{1}{\kappa}}+A_2^{\frac{1}{\kappa}})^{\kappa}}{\kappa}
\epsilon^{\kappa}
+o(\epsilon^{\kappa})\Label{4-5-2}
\end{align}
hold.

When the function $x \mapsto \log f(x)$ is concave,
the relations (\ref{4-5-3}) and (\ref{4-5-4}) 
yield that
\begin{align}
\alpha(\vec{\theta}_{ML},\theta)\ge
\frac{\overline{\alpha}_1(\theta)}{2^{\kappa}}.
\end{align}

Further, when $f(x)$ is monotonically decreasing, 
Lemma \ref{4-5-7} and (\ref{4-5-6}) imply that
\begin{align}
\beta(\vec{\theta}_{ML},\theta,\epsilon)
\cong\frac{A_1}{\kappa_1}\epsilon^{\kappa_1}+ o(\epsilon^{\kappa_1}).
\end{align}

\subsection{Exponential rates of estimators for interval estimations}
In order to improve the exponential rate for a fixed width $2\epsilon$,
we focus on 
the likelihood ratio estimator, which is discussed 
by Huber, Sievers, and Fu 
from the large deviation viewpoint 
concerning the regular family\cite{Hu,Sie,Fu85}.
Assume that the function $d_{\epsilon}(x):=f(x + \epsilon) /f(x -\epsilon)$ is monotonically decreasing w.r.t. $x$ when both $f(x + \epsilon)$ and $f(x -\epsilon)$ are not zero.
Then, we can define 
the likelihood ratio estimator 
$\vec{\theta}_{LR,\epsilon}:
=\{\theta_{LR,\epsilon,n}(x_1 , \ldots , x_n)\}$,
which depends on the constant $\epsilon \,> 0$, as shown by
\begin{align}
\theta_{LR,\epsilon,n}:=
\frac{1}{2} \left(\sup\{ z | k(z) \,< 0 \}+
\inf \{ z | k(z) \,> 0 \} \right), \Label{725.2}
\end{align}
where the monotonically decreasing function $k(z)$ is defined by
\begin{align}
k(z):= \frac{1}{n} \sum_{i=1}^n \left(
\log f( x_i - z +\epsilon ) - 
\log f( x_i - z -\epsilon )\right)
\Label{18.1}.
\end{align}
Note that when $\log f(x)$ is concave, the above condition is satisfied. This definition is well defined although the monotonically decreasing function $k(z)$ is not continuous.

If the support of $f$ is $(a,b)$, we need to modify the definition as follows. 
In this case, we modify the estimator $\theta_{LR,\epsilon,n}$ 
by using the two estimators 
$\overline{\theta}_n$ and $\underline{\theta}_n$. 
When $\underline{\theta}_n- \overline{\theta}_n \,> 2 \epsilon$, 
the estimated value is defined by (\ref{725.2}) in the interval $(\underline{\theta}_n- \epsilon ,\overline{\theta}_n + \epsilon)$. 
When $\underline{\theta}_n- \overline{\theta}_n \le 2 \epsilon$, 
we define $\theta_{LR,\epsilon,n}:= 
\frac{1}{2}\left(\underline{\theta}_n+ \overline{\theta}_n\right)$. 
Moreover, when the support of $f$ is a half line $(0, \infty)$, 
the estimated value is defined by (\ref{725.2}) in the half line $(\underline{\theta}_n- \epsilon , \infty)$.

Then, the exponential rate is characterized as follows.
(A regular version of this lemma was discussed by Huber \cite{Hu}, 
Sievers \cite{Sie}, and Fu \cite{Fu85}.)
\begin{lem}\Label{thm7.2}
When $\log f(x)$ is concave, the equation 
\begin{align}
\min\{ 
\beta^-( \vec{\theta}_{LR,\epsilon}, 
\theta , \epsilon ),
\beta^+( \vec{\theta}_{LR,\epsilon}, 
\theta , \epsilon )\}
= \sup_{0\,< s \,< 1}
I^s(f_{\theta-\epsilon }\| f_{\theta+\epsilon}) . 
\Label{725.3}
\end{align}
holds, where $f_{\theta}(x):= f(x- \theta)$.
\end{lem}
Therefore, the equality of inequality (\ref{7-a}) holds in this case.

\begin{proof}
Note that 
$\beta^+( \vec{\theta}_{LR,\epsilon}, 
\theta , \epsilon )
=\beta^+( \vec{\theta}_{LR,\epsilon}, 
\theta - \epsilon , 
\epsilon)$
and 
$\beta^-( \vec{\theta}_{LR,\epsilon},
 \theta , \epsilon )
=\beta^-( \vec{\theta}_{LR,\epsilon}, 
\theta + \epsilon , 
\epsilon )$
because of the shift-invariance. From the concavity, the condition 
$\theta_{LR,\epsilon,n} \le \theta$ is equivalent to the condition 
$\sup\{ z | k(z) \,< 0\} \le \theta$, which implies that $k(\theta)\ge 0$.
Thus, we have 
\begin{align*}
\sum_i
\log f(x_i - \theta +\epsilon)
- \log f(x_i - \theta -\epsilon)
\ge 0.
\end{align*}
Conversely,
the condition \begin{align*}
\sum_i
\log f(x_i - \theta +\epsilon)
- \log f(x_i - \theta -\epsilon)
\,> 0
\end{align*}
implies that
$k(\theta) \,> 0$.
Thus, we have
$\sup\{ z | k(z) \,< 0\} \le \theta$,
which is equivalent to the condition $\theta_{LR,\epsilon,n} \le \theta$. 
Therefore, we have the relations 
\begin{align*}
\{ f_{\theta-\epsilon}^n (\vec{x}_n) \,>
f_{\theta+\epsilon}^n (\vec{x}_n) \}
\subset 
\{ \theta_{LR,\epsilon,n} \le \theta\}
\subset
\{ f_{\theta-\epsilon}^n (\vec{x}_n) \ge
f_{\theta+\epsilon}^n (\vec{x}_n) \}. 
\end{align*}
Similarly, we can prove
\begin{align*}
\{ f_{\theta+\epsilon}^n (\vec{x}_n) \,>
f_{\theta-\epsilon}^n (\vec{x}_n) \}
\subset 
\{ \theta_{LR,\epsilon,n} \ge \theta\}
\subset
\{ f_{\theta+\epsilon}^n (\vec{x}_n) \ge
f_{\theta-\epsilon}^n (\vec{x}_n) \} .
\end{align*}
Applying (\ref{11}), we obtain
\begin{align*}
\min\{ 
\beta^-( \vec{\theta}_{LR,\epsilon}, 
\theta + \epsilon, 
\epsilon ),
\beta^+( \vec{\theta}_{LR,\epsilon}, 
\theta - \epsilon,
 \epsilon )\}
= \sup_{0\,< s \,< 1}
I^s(f_{\theta-\epsilon }\| f_{\theta+\epsilon}),
\end{align*}
which implies equation (\ref{725.3}).
\end{proof}

Next, 
in order to improve the estimator $\vec{\underline{\theta}}$ for a fixed 
value $\epsilon$,
we define
the estimator 
$\vec{\underline{\theta}}_\epsilon 
:= \{ \underline{\theta}_{\epsilon,n} := \underline{\theta}_n - \epsilon\}$.
This estimator satisfies the following lemma.
\begin{lem}
The estimator $\vec{\underline{\theta}}_\epsilon$ satisfies 
\begin{align}
\beta^+(\vec{\underline{\theta}}_\epsilon, 
\theta, \epsilon)
&=
- \log \left(\int_{a}^{b-2 \epsilon}f(x) \,d x\right).\Label{3-8-1-a}\\
\beta^-(\vec{\underline{\theta}}_\epsilon, \theta, \epsilon)
&=\infty \Label{3-8-1-b}.
\end{align}
Further, when $f(x)$ is monotonically decreasing, 
the 
exponential rate of 
the estimator $\vec{\underline{\theta}}_\epsilon$ 
has another form
\begin{align}
\beta^+(\vec{\underline{\theta}}_\epsilon, 
\theta, \epsilon)
&=\sup_{0\,< s \,< 1} I^s (f_{\theta-\epsilon}
\|f_{\theta+\epsilon}) \Label{3-8-1}.
\end{align}
\end{lem}
\begin{proof}
From the construction of $\vec{\underline{\theta}}_\epsilon$,
the relations (\ref{3-8-1-a}) and (\ref{3-8-1-b}) follow
from (\ref{3-8-14}).
Then, we proceed to the proof of (\ref{3-8-1}).
Since $\underline{\theta}_n \,> \theta$, we have $\underline{\theta}_{\epsilon,n} \,> \theta- \epsilon$, which implies (\ref{3-8-1}). If $\underline{\theta}_{\epsilon,n} \ge \theta +\epsilon$, we have $\underline{\theta}_{n} \ge \theta + 2 \epsilon$.
Thus, $f(x_i -(\theta +2\epsilon)) 
\ge f(x_i -\theta )$ for any $i= 1, \ldots, n$.
Therefore, 
\begin{align*}
f^n_{\theta+2\epsilon} (\vec{x}_n)
\ge f^n_{\theta} (\vec{x}_n).
\end{align*}
Conversely, if $f^n_{\theta+2\epsilon} (\vec{x}_n) \ge f^n_{\theta} (\vec{x}_n)$,
we have $\underline{\theta}_{n} \ge \theta + 2 \epsilon$. Thus, 
\begin{align*}
f^n_{\theta+2\epsilon} \{ \vec{x}_n|
f^n_{\theta+2\epsilon} (\vec{x}_n)
\,< f^n_{\theta} (\vec{x}_n)\}
= 
f^n_{\theta+2\epsilon} \{ \vec{x}_n|
\underline{\theta}_{n} 
\,< \theta + 2 \epsilon\}=0.
\end{align*}
Since the likelihood test $\{ \vec{x}_n| \underline{\theta}_{n} \ge \theta + 2 \epsilon\}$ achieves the optimal rate (\ref{11}), we have 
\begin{align*}
&\lim -\frac{1}{n}\log
f^n_{\theta} \{ \vec{x}_n|
f^n_{\theta+2\epsilon} (\vec{x}_n)
\ge f^n_{\theta} (\vec{x}_n)\} \\
=&
\lim -\frac{1}{n}\log
\left(
f^n_{\theta} \{ \vec{x}_n|
f^n_{\theta+2\epsilon} (\vec{x}_n)
\ge f^n_{\theta} (\vec{x}_n)\}+
f^n_{\theta+2\epsilon} \{ \vec{x}_n|
f^n_{\theta+2\epsilon} (\vec{x}_n)
\,< f^n_{\theta} (\vec{x}_n)\}\right)\\
=&
\sup_{0\,< s \,<1}
I^s(f_{\theta}\|f_{\theta+2\epsilon}) 
= 
\sup_{0\,< s \,<1}
I^s(f_{\theta- \epsilon}\|f_{\theta+\epsilon}) .
\end{align*}
\end{proof}

\subsection{Slopes of exponential rates of estimators useful for interval estimation}
Next, we proceed to the slopes of exponential rates.
Concerning
the estimator 
$\vec{\underline{\theta}}_{\epsilon}$,
from (\ref{3-8-1-a}) and (\ref{3-8-1-b}),
we obtain
the following characterization:
\begin{align}
\beta(\vec{\underline{\theta}}_{\epsilon}, \theta, \epsilon )
\cong
A_1 \frac{2^{\kappa_1}}{\kappa_1}
\epsilon^{\kappa_1}+o(\epsilon^{\kappa_1}).\Label{4-5-9}
\end{align}
Further, when $f(x)$ is monotonically decreasing, 
the equation (\ref{3-8-1}) implies
the equation 
\begin{align}
\lim_{\epsilon\to +0}
\frac{\beta(\vec{\underline{\theta}}_{\epsilon}, \theta, \epsilon )}
{g(\epsilon)}
= \overline{\alpha}_1(\theta).
\end{align}
When the function $\log f(x)$ is concave, 
Lemma \ref{thm7.2} yields 
the equation 
\begin{align}
\lim_{\epsilon\to +0}
\frac{\beta(\vec{\theta}_{LR,\epsilon}, \theta, \epsilon )}
{g(\epsilon)}
= \overline{\alpha}_1(\theta).
\end{align}
Thus, the following theorem holds.
\begin{thm}
When the function $\log f(x)$ is concave or monotonically decreasing,
the relation 
\begin{align*}
\alpha_1(\theta)= \overline{\alpha}_1(\theta).
\end{align*}
holds.
\end{thm}
Therefore, 
if the condition of the above theorem holds and 
condition (\ref{16.3}) in Theorem \ref{thm3} does not hold, i.e., 
$I^s_{g,\theta}$ is not symmetric,
then the first criterion is different from the second one, i.e.,  
\begin{align}
\alpha_1(\theta)> \alpha_2(\theta).
\end{align}

\section{Conclusion}
We have discussed large deviation theories under a more general setting
than existing studies. 
For this purpose, we introduced two criteria for the bound of estimation 
accuracy from the large deviation viewpoint.
One criterion $\alpha_1(\theta)$
corresponds to the interval estimation with taking 
the limit that the width of error bar goes to $0$.
The other $\alpha_2(\theta)$
corresponds to the point estimation.
The upper bounds of them are given by the limit of the relative R\'{e}nyi entropy.
These characterizations have been derived by the method of 
simple hypothesis testing.
We have also calculated the slope of the exponential decreasing rates of 
several estimators.
As a result, 
we have succeeded in deriving upper bounds $\overline{\alpha}_1(\theta)$ and $\overline{\alpha}_2(\theta)$
of $\alpha_1(\theta)$ and $\alpha_2(\theta)$ and
a necessary and sufficient condition for the coincidence of these two upper bounds $\overline{\alpha}_1(\theta)$ and $\overline{\alpha}_2(\theta)$. 
In the next step, we have treated several estimators as candidates to attain the optimal bounds
$\alpha_1(\theta)$ and $\alpha_2(\theta)$ in a local shift family.
That is, we derived lower bounds of $\alpha_1(\theta)$ and $\alpha_2(\theta)$.
Further, we derived a sufficient condition for 
the gap between two criteria $\alpha_1(\theta)$ and $\alpha_2(\theta)$.

Unfortunately, we cannot calculate 
two criteria $\alpha_1(\theta)$ and $\alpha_2(\theta)$
in the concrete examples.
For this purpose, we need to calculate 
the limit of the relative R\'{e}nyi entropy,
which was discussed in another paper \cite{Haya-2}.
In the next paper, we will treat this calculation based on the obtained result and the result by Hayashi \cite{Haya-2}.

Historically, Nagaoka initiated a discussion of two kinds 
of large deviation bounds in a quantum setting \cite{Naga2,Naga3}, 
and Hayashi discussed these in more depth \cite{Haya}. 
The two kinds of large deviation bounds do not necessarily coincide 
in a quantum setting. 
In quantum setting, this difference corresponds to the non-uniqueness of quantum 
extension of Fisher information.
In the quantum case, $\alpha_1(\theta)$ corresponds to 
Kubo-Mori-Bogoljubov (KMB) inner product and $\alpha_2(\theta)$ does to Symmetric logarithmic derivative (SLD) inner product.
This research is strongly motivated by this quantum study.
In the quantum case, 
the family of estimators attaining 
the bound $\overline{\alpha}_1(\theta)$
depends on the true parameter.
However, 
in some of non-regular location families, 
such a family of estimators does not depends on the true parameter.
This is different point between our setting and quantum setting.
Gaining an understanding of 
these differences from a unified viewpoint remains a goal for the future.

\section*{Acknowledgment}
The present study was supported in part by 
Laboratory for Mathematical Neuroscience, 
Brain Science Institute, RIKEN
and MEXT through a Grant-in-Aid for Scientific Research on Priority Area ``Deepening and Expansion of Statistical Mechanical Informatics (DEX-SMI),'' No. 18079014. 
The author is grateful for Professor Shun-ichi Amari
to helpful discussions on this topics.

\appendix

\section{Concave function}
\begin{lem}\Label{concave}
When a concave function $f \ge 0$ is defined in $(0,1)$,
\begin{align*}
\inf_{x \ge0}
sx +(1-s)\sup_{0 \,< t \,< 1}
\frac{-tx +f(t)}{1-t}
=
\inf_{x \,>0}\sup_{0 \,< t \,< 1}
\frac{(s-t)x +(1-s)f(t)}{1-t}
= f(s).
\end{align*}
\end{lem}
\begin{proof}
Substituting $s$ into $t$ we have
\begin{align*}
f(s) \le
\sup_{0 \,< t \,< 1}
\frac{(s-t)x +(1-s)f(t)}{1-t}.
\end{align*}
Taking the infimum $\inf_{x \,>0}$, we obtain 
\begin{align*}
f(s) \le \inf_{x \,>0}\sup_{0 \,< t \,< 1}
\frac{(s-t)x +(1-s)f(t)}{1-t}.
\end{align*}

Next, we proceed to the opposite inequality. From the concavity of $f$, we can define the upper derivative $\overline{f'}$ and the lower derivative $\underline{f'}$ as
\begin{align*}
\overline{f'}(s):= \lim_{\epsilon to +0}
\frac{f(s)- f(s-\epsilon)}
{\epsilon} ,\quad 
\underline{f'}(s):= \lim_{\epsilon to +0}
\frac{f(s+\epsilon)- f(s)}
{\epsilon} 
\end{align*}
Since the concavity guarantees that
\begin{align*}
\frac{f(s+\epsilon)- f(s)}
{\epsilon} 
\le \underline{f'}(s) \le \overline{f'}(s)\le 
\frac{f(s)- f(s-\epsilon)}
{\epsilon} , \quad \forall \epsilon\,> 0 ,
\end{align*}
we obtain
\begin{align*}
&\frac{(s-s)(f(s)+(1-s)\overline{f'}(s))
+(1-s)f(s)}{1-s}
- 
\frac{(s-t)(f(s)+(1-s)\overline{f'}(s))
+(1-s)f(t)}{1-t} \\
=& \frac{1-s}{1-t}
\left( f(s)-f(t) +(t-s)\overline{f'}(s)\right)
\ge 0,\quad \forall s, \forall t \in (0,1).
\end{align*}
Therefore,
\begin{align*}
f(s) \ge &\sup_{0\,< t \,< 1}
\frac{(s-t)(f(s)+(1-s)\overline{f'}(s))
+(1-s)f(t)}{1-t} \\
\ge &
\inf_{x \,> 0}
\sup_{0\,< t \,< 1}
\frac{(s-t)(f(s)+(1-s)\overline{f'}(s))
+(1-s)f(t)}{1-t}.
\end{align*}
The proof is now complete.
\end{proof}

\section{Other lemmas}\Label{A}
\begin{lem}\Label{ap1}
When $g$ is strictly monotonically decreasing and continuous,
$g(0)=0$, and the limit 
$\lim_{\epsilon \to +0}\frac{g(x\epsilon)}{g(\epsilon)}$ exists, 
there exists $\kappa \,> 0$ such that
\begin{align}
x^\kappa = 
x \,> 0. \Label{82}
\end{align}
\end{lem}
\begin{proof}
Let $h(x)$ be the RHS of (\ref{82}).
Since 
\begin{align*}
\lim_{\epsilon \to +0}
\frac{g(xy \epsilon)}{g(\epsilon)}
=
\lim_{\epsilon \to +0}
\frac{g(xy \epsilon)}{g(y\epsilon)}
\lim_{\epsilon \to +0}
\frac{g(y \epsilon)}{g(\epsilon)},
\end{align*}
$h(xy)=h(x)h(y)$.
Thus, there exists $\kappa \,> 0$ satisfying (\ref{82}).
\end{proof}

\begin{lem}\Label{L38}
For $0 \,< s \,< \frac{1}{2}$, the inequalities
\begin{align*}
2s I^{\frac{1}{2}}(p\|q)
\le I^{s}(p\|q)
\le 2(1-s) I^{\frac{1}{2}}(p\|q)
\end{align*}
hold.
\end{lem}
\begin{proof}
Since $\left(\frac{1}{2s}\right)^{-1}+
\left(\frac{1}{1-2s}\right)^{-1}=1$,
the H\"{o}lder inequality guarantees that
\begin{align*}
\int_{\Omega}
p^{s}(\omega) q^{1-s}(\omega) \,d \omega 
=&
\int_{\Omega}
\left(p^{s}(\omega) q^{s}(\omega) \right)
\left(q^{1-2s}(\omega)\right) \,d \omega \\
\le &
\left(\int_{\Omega}
\left(p^{s}(\omega) q^{s}(\omega) \right)^{\frac{1}{2s}}
\,d \omega \right)^{2s}\cdot
\left(
\int_{\Omega}
\left(q^{1-2s}(\omega)\right)^{\frac{1}{1-2s}}
 \,d \omega 
\right)^{1-2s}\\
=& \left(\int_{\Omega}
p^{\frac{1}{2}}(\omega) q^{\frac{1}{2}}(\omega) 
\,d \omega \right)^{2s}.
\end{align*}
Thus, we obtain 
\begin{align*}
I^{s}(p\|q) \ge 2s I^{\frac{1}{2}}(p\|q).
\end{align*}
Similarly, since 
$\left( 2(1-s)\right)^{-1}+
\left( \frac{2(1-s)}{1-2s}\right)^{-1}=1$,
we can apply the H\"{o}lder inequality as
\begin{align*}
&\int_{\Omega}
p^{\frac{1}{2}}(\omega) q^{\frac{1}{2}}(\omega) 
\,d \omega \\
=&\int_{\Omega}
\left(p^{s}(\omega) q^{1-s}(\omega) 
\right)^{\frac{1}{2(1-s)}}
p^{\frac{1-2s}{2(1-s)}}(\omega) 
\,d \omega \\
\le &
\left(\int_{\Omega}
\left(p^{s}(\omega) q^{1-s}(\omega) 
\right)^{\frac{1}{2(1-s)}\cdot 2(1-s)}
\,d \omega \right)^{\frac{1}{2(1-s)}}\cdot
\left(\int_{\Omega}
p^{\frac{1-2s}{2(1-s)}\cdot\frac{2(1-s)}{1-2s}}(\omega) 
\,d \omega \right)^{\frac{1-2s}{2(1-s)}} \\
=&
\left(\int_{\Omega}
p^{s}(\omega) q^{1-s}(\omega) \,d \omega
\right)^{\frac{1}{2(1-s)}},
\end{align*}
which implies that
\begin{align*}
I^{\frac{1}{2}}(p\|q) \ge \frac{1}{2(1-s)}
I^{s}(p\|q).
\end{align*}
\end{proof}

\end{document}